\documentclass[12pt]{article}
\usepackage{amssymb,amsmath,amstext}
\usepackage{a4}
\usepackage{colordvi}
\newtheorem{corollary}{Corollary}[section]
\newtheorem{step}{Claim}[section]
\newtheorem{definition}{Definition}[section]
\newtheorem{lemma}{Lemma}[section]

\newcommand{\qed}{\hbox to .60em{\vrule width .60em height .60em}}

\newtheorem{proposition}{Proposition}[section]
\newtheorem{remark}{Remark}[section]
\newtheorem{theorem}{Theorem}[section]

\newcommand{\ga}{\gamma}
\newcommand{\ze}{\zeta}

\newcommand{\al}{\alpha}
\newcommand{\be}{\beta}
\newcommand{\na}{\nabla}

\newcommand{\bstep}{\begin{step}}
\newcommand{\estep}{\end{step}}
\newcommand{\blem}{\begin{lemma}}
\newcommand{\elem}{\end{lemma}}
\newcommand{\brem}{\begin{remark}}
\newcommand{\erem}{\end{remark}}
\newcommand{\bthm}{\begin{theorem}}
\newcommand{\ethm}{\end{theorem}}
\newcommand{\beqn}{\begin{equation}}
\newcommand{\eeqn}{\end{equation}}
\newcommand{\eeq}{\end{equation}}
\newcommand{\beq}{\begin{equation}}
\newcommand{\eitem}{\end{itemize}}
\newcommand{\bitem}{\begin{itemize}}
\newcommand{\eenum}{\end{enumerate}}
\newcommand{\benum}{\begin{enumerate}}

\newcommand{\ds}{\displaystyle}

\def\E{\mathbb{E}}
\def\F{\mathbb{F}}

\def\MM{\mathbb{M}}

\def\O{\mathbb{O}}
\def\P{\mathbb{P}}
\def\R{\mathbb{R}}
\def\S{\mathbb{S}}

\def\C{\mathbb{C}}
\def\G{\mathbb{G}}

\def\trace{\mbox{Trace}}

\def\a{\alpha}
\def\b{\beta}

\def\T{\mathcal{T}}

\title{Newton's Method on Riemannian Manifolds:\\
Covariant Alpha-Theory.}

\author{Jean-Pierre Dedieu%
\thanks{%
MIP. D\'epartement de Math\'ematique,
Universit\'e Paul Sabatier,
31062 Toulouse cedex 04, France
({\tt dedieu@mip.ups-tlse.fr}).
}
\and Pierre Priouret%
\thanks{%
MIP. D\'epartement de Math\'ematique,
Universit\'e Paul Sabatier,
31062 Toulouse cedex 04, France
({\tt priouret@mip.ups-tlse.fr}).
}
\and Gregorio Malajovich%
\thanks{%
Departamento de Matem\'atica Aplicada,
Universidade Federal de Rio de Janeiro,
Caixa Postal 68530,
CEP 21945-970, Rio de Janeiro, RJ, Brazil
({\tt gregorio@labma.ufrj.br}).
}}

\date{January 15, 2003}

\begin{document}
\maketitle
\begin{abstract}
In this paper, Smale's
$\alpha$ theory is generalized to the context of intrinsic
Newton iteration on geodesically complete analytic 
Riemannian and Hermitian manifolds. Results are valid for
analytic mappings from a manifold to a linear space of the same
dimension, or for analytic vector fields on the manifold.
The invariant $\gamma$ is defined by means of high order
covariant derivatives. Bounds on the size of the basin
of quadratic convergence are given. If the ambient manifold
has negative sectional curvature, those bounds depend on
the curvature. A criterion of quadratic convergence for Newton
iteration from the information available at a point is also given.
\end{abstract}

\section{Introduction and main results.}

Numerical problems posed in manifolds arise in many natural contexts. Classical examples 
are given by the eigenvalue problem, the symmetric eigenvalue problem, invariant subspace computations, 
minimization problems with orthogonality constraints, optimization problems with 
equality constraints ... etc. In the first example, $Ax = \lambda x$, the unknowns are the eigenvalue $\lambda \in \C$ and the eigenvector $x \in \P_{n-1}(\C)$, the complex projective space consisting of complex vector lines through
the origin in $\C^n$. In the second example, $Ax = \lambda x$, $A$ real and symmetric, the  unknowns are $\lambda \in \R$ and $x \in \S^{n-1}$, the unit sphere in $\R^n$. In the third example the unknown is a $k-$dimensional subspace 
contained in $\C^n$ that is an element of the Grassmann manifold $\G_{n,k}(\C )$. The fourth example involves the orthogonal group, the special orthogonal group or the Stiefel manifold ($n \times k$ matrices with orthonormal columns). The last example leads to problems posed on submanifolds in $\R^n$.

For such or similar problems our objective is to design algorithms which respect their geometrical structure. 
We follow here the lines of the Geometric Integration Interest Group (http://www.focm.net/gi/) who showed the interest of such an approach.

The first author's original motivation came from homogeneous and multihomogeneous polynomial systems (Dedieu-Shub \cite{Dedieu-Shub-2}) and also from a model for the human spine (Adler-Dedieu-Margulies-Martens-Shub~\cite{ADMMS}) with configuration space $SO(3)^{18}$. A second motivation, for the second author, came from sparse polynomial systems of equations where the solutions belong to a certain toric variety (Malajovich-Rojas \cite{MR}).

For such problems one often has to compute the solutions of a system of equations or to find the zeros of a vector field. For this reason we investigate here one of the most famous method to approximately solve these problems: the Newton method.

In this paper, we investigate the local behavior of Newton's iteration
close to a solution. While a lot is known about Newton's iteration in
linear spaces~\cite{BCSS}, little is known about intrinsic Newton's
iteration in more general manifolds. 
Our main results here (Theorems~\ref{th3} to~\ref{thV4} 
below) extend Smale's $\alpha$-theory to analytic Riemannian manifolds.
$\alpha$ theory provides a criterion for the quadratic convergence of
Newton's iteration in a neighborhood of a solution. This criterion
depends on available data at the approximate solution. One important
application (out of the scope of this paper) is the construction of
rigorous homotopy algorithms for the solution of non-linear equations.

More precisely, we will study quantitative aspects of Newton's method for
finding zeros of mappings $f: \MM_n \rightarrow \R^n$ and vector fields
$X: \MM_n \rightarrow T\MM_n$. Here $\MM_n$ denotes a real complete analytic 
Riemannian manifold, $T\MM_n$ its tangent bundle, $f$ and $X$ are analytic. 
We denote by $T_z\MM_n$ the tangent space at $z$ to $\MM_n$, by 
$\langle .,. \rangle_z$ the scalar product on $T_z\MM_n$ with associated norm
$\| . \|_z$, by $d$ the Riemannian metric on $\MM_n$ and by 
$\exp_z : T_z\MM_n \rightarrow \MM_n$ 
the exponential map. This map is defined on the whole tangent bundle $T\MM_n$ 
because $\MM_n$ is assumed to be complete. We denote by ${\bf r}_z > 0$ the
radius of injectivity of the exponential map at $z$. Thus, 
$\exp_z : B_{T_z}(0,{\bf r}_z) \rightarrow B_{\MM_n}(z,{\bf r}_z)$ is one to one ($B(u,r)$ is 
the open ball about $u$ with radius $r$, $\bar B(u,r)$ is the closed ball). 

When $\MM_n = \R^n$ the Newton operator associated with $f$ is defined by 
$$N_f(z) = z - Df(z)^{-1}f(z).$$ 
In this context $T_z {\R^n}$ may be identified to $\R^n$
and $\exp_z (u) = z + u$ so that 
$$N_f(z) = \exp_z (-Df(z)^{-1}f(z)).$$ 
This formula makes sense in the context of Riemannian manifolds and we define 
the Newton operator $N_f : \MM_n \rightarrow \MM_n$ in this way. 

When, instead of a mapping $\MM_n \rightarrow \R^n$ we consider a vector field
$X: \MM_n \rightarrow T\MM_n$, in order to define Newton's method, we resort to
an object studied in differential geometry; namely, the covariant derivative
of vector fields. Let $\nabla$ denote the Levi-Civita connection on $\MM_n$. 
For any vector fields $X$ and $Y$ on $\MM_n$, $\na_X(Y)$ is called the covariant 
derivative of $Y$ with respect to $X$. Since $\nabla$ is tensorial in $X$
the value of $\na_X(Y)$ at $z \in \MM_n$ depends only on the tangent vector
$u = X(z) \in T_z\MM_n$. For this reason we denote it
$$(\na_X(Y))(z) = DY(z)(u).$$
It is a linear map 
$$DY(z) : T_z\MM_n \rightarrow T_z\MM_n.$$
The Newton operator for the vector field $X$ is defined by
$$N_X(z) = \exp_z (-DX(z)^{-1}X(z)).$$ 
Notice this definition coincides with the usual one when $X$ is a vector
field in $\R^n$ because the covariant derivative is just the usual derivative.

In a vector space framework, Newton's method makes zeros of $f$ with non-singular 
derivative correspond to fixed points of $N_f$ and 
Newton sequences $x_{k+1} = N_f(x_k)$, for an initial point $x_0$ taken close
to such a fixed point $\ze$, converge quadratically to $\ze$. In this paper, 
our aim is to make these statements precise in our new geometric framework and 
to investigate quantitative aspects.
We have in mind the following two theorems which are valid when $\MM_n$ is equal
to $\R^n$ or in the more general context of an analytic mapping $f : \E \rightarrow \F$ between two real or complex
Banach spaces:

\begin{theorem} \label{th1} {\bf ($\ga-$Theorem, Smale, 1986)} Suppose that $f(\ze) = 0$ and
$Df(\ze)$ is an isomorphism. Let
$$\ga(f,z) = \sup_{k \geq 2} \left \| Df(z)^{-1} \frac{D^kf(z)}{k!} \right 
\|^{1/k-1}.$$
If 
$$\| z - \ze \| \leq \frac{3 - \sqrt{7}}{2 \ga(f,\ze)}$$
then the Newton sequence $z_k = N_f^{(k)} (z)$ is defined for all $k \geq 0$
and 
$$\| z_k - \ze \| \leq \left ( \frac{1}{2} \right )^{2^k - 1}
\| z - \ze \|.$$
\end{theorem}

For a proof see Blum-Cucker-Shub-Smale \cite{BCSS} Chap. 8, Theorem 1. 
The second theorem we want to extend to the context of Riemannian manifolds
is the following:

\begin{theorem} \label{th2} {\bf ($\a-$Theorem, Smale, 1986)} Let 
$$\be(f,z) = \| Df(z)^{-1}f(z) \|$$ and
$$\a(f,z) = \be(f,z) \ga(f,z).$$ 
We also let $\a(f,z) = \infty$ when $Df(z)$ is not invertible. There is a universal 
constant $\a_0 > 0$ with the following property: if $\a(f,z) < \a_0$ then 
there is a zero $\ze$ of $f$ such that $Df(\ze)$ is an isomorphism
and such that the Newton sequence $z_k = N_f^{(k)} (z)$ is defined for all 
$k \geq 0$ and satisfies
$$\| z_k - \ze \| \leq \left ( \frac{1}{2} \right )^{2^k - 1}
\| z - \ze \|.$$
Moreover, the distance from $z$ to the zero $\ze$ is at most $2\be(f,z)$.
\end{theorem}

This second theorem is proved in Smale \cite{Smale86} with the constant 
$\a_0 = 0.13071 \ldots$ and Kim \cite{Kim85} and \cite{Kim88} for a one-dimensional 
version.

\subsection{Definitions and notations.}

In order to generalize these two results we have to define the corresponding invariants in the context of Riemannian manifolds. The material contained in this section is classical in Riemannian geometry. 
The reader is refered to a textbook on this subject, for example: Dieudonn\'e
\cite{Dieud}, Do Carmo \cite{Docarmo}, Gallot-Hulin-Lafontaine \cite{GHL}, 
Helgason \cite{Helgason}, O'Neill \cite{Oneill}.

\begin{definition} {\bf (Tensors.)} The space of $p-$contravariant and $q-$covariant {analytic tensor fields} 
$$T : T(\MM_n)^p \times T^*(\MM_n)^q \rightarrow {\cal F}(\MM_n)$$
is denoted by $\T_q^p(\MM_n)$. An $m-$tuple of such tensor fields is called a vectorial tensor field and 
the space of vectorial tensor fields is denoted by $\T_q^p(\MM_n, \R^m)$.
\end{definition}

Here $T^*(\MM_n)$ is the cotangent bundle on $\MM_n$ (the space of $1-$forms) and ${\cal F}(\MM_n)$ 
the space of scalar analytic functions defined on $\MM_n$. We let ${\cal F}(\MM_n) = \T_0^0(\MM_n)$.
Let $\nabla$ denote the Levi-Civita connection on $\MM_n$. For any vector field
$X$ and $Y$ on $\MM_n$, $\na_X(Y)$ is called the covariant derivative of $Y$ with
respect to $X$. 

\begin{definition} {\bf (Covariant derivative for tensor fields.)} Let $X$ be a vector field on $\MM_n$. For any integers $p$, $q \ge 0$ and any tensor field $T \in \T_q^p(\MM_n)$ the covariant derivative is defined by:
\begin{itemize}\item $\na_X(g) = X(g) = Dg (X)$ the derivative of $g$ 
along the vector field $X$ when $g$ is a function: $g \in \T_0^0(\MM_n)$
\item $\na_X(Y)$ is given by the connection when $Y$ is a vector field i.e.
$Y \in \T_0^1(\MM_n)$
\item For a $1-$form $\omega \in \T_1^0(\MM_n)$ its covariant derivative is 
the $1-$form defined by 
$$\na_X(\omega)(Y) = X(\omega (Y)) - \omega (\na_X(Y))$$
for any vector field $Y$.
\item For a tensor field $T \in \T_q^p(\MM_n)$ the covariant derivative is 
the tensor field $\na_XT \in \T_q^p(\MM_n)$ defined by 
$$\na_XT(\omega^1 \ldots \omega^p,Y_1 \ldots Y_q) = 
X(T(\omega^1 \ldots \omega^p,Y_1 \ldots Y_q)) - $$
$$T(\na_X(\omega^1) \ldots \omega^p,Y_1 \ldots Y_q) - \ldots -
T(\omega^1 \ldots \omega^p,Y_1 \ldots \na_X(Y_q))$$
for any $1-$forms $\omega^i$ and vector fields $Y_j$.
\item For a vectorial tensor field $T \in \T_q^p(\MM_n, \R^m)$
$$\na_X \left ( \displaystyle
\begin{array}{cc}
T_1\\
\vdots\\
T_m\\
\end{array} \right ) = 
\left ( \displaystyle
\begin{array}{cc}
\na_X T_1\\
\vdots\\
\na_X T_m\\
\end{array} \right ).$$
\end{itemize}
\end{definition}

\begin{definition} {\bf (Covariant $k-$th derivative for tensor fields.)} Let $X$ be a vector field on $\MM_n$. For any integers $p$, $q \ge 0$ and any tensor fields $T \in \T_q^p(\MM_n, \R^m)$ the $k-$th covariant derivative is defined 
inductively by 
$$\na_X^k T = \na_X \left(\na_X^{k-1} T \right).$$
\end{definition}

Since the covariant derivative is tensorial in $X$, its value at a given point $z \in \MM_n$ depends only on the 
vector $X(z)$. For this reason, the following definition makes sense:

\begin{definition} {\bf (Covariant $k-$th derivative for tensor fields at a point.)} \label{dkf} Let a point $z \in \MM_n$ and a vector $u \in \T_z(\MM_n)$ be given. Let $X$ be a vector field such that $X(z) = u$.
For any integers $p$, $q \ge 0$ and any tensor field $T \in \T_q^p(\MM_n, \R^m)$ the value at $z$ of the $k-$th covariant derivative is denoted by:
$$D^kT(z)(u, \ldots ,u) = D^kT(z)u^k = (\na_X^k T) (z).$$
It defines a $k-$multilinear map
$$D^kT(z) : \left( \left( T_z\MM_n \right)^p \times \left( T^*_z\MM_n \right)^q \right)^k \rightarrow \R^m.$$
\end{definition}

\begin{definition} {\bf (Norm of a multilinear map.)} \label{norm} Let 
$$M : \left(T_z\MM_n\right)^k \rightarrow \R^m$$
be a $k-$multilinear map.
Its norm is defined by
$$\| M \|_z = \sup \| M(u_1, \ldots , u_k) \|_{\R^m}$$
where the supremum is taken for all the vectors $u_j \in T_{z}\MM_n$ such that $\| u_j \|_z = 1.$
\end{definition}

The following definition extends the definition of $\gamma (f,z)$ to a Riemannian context.

\begin{definition} {\bf (Gamma.)} Let a map $f : \MM_n \rightarrow \R^n$ and a vector field $X : \MM_n \rightarrow T\MM_n$ be given. For any point $z \in \MM_n$ we let
$$\ga(f,z) = \sup_{k \geq 2} \left \| Df(z)^{-1} \frac{D^kf(z)}{k!} \right \|_z^{1/k-1},$$
$$\ga(X,z) = \sup_{k \geq 2} \left \| DX(z)^{-1} \frac{D^kX(z)}{k!} \right \|_z^{1/k-1}.$$
We also let $\ga(f,z) = \infty$ when $Df(z)$ is not invertible, idem for $\ga(X,z)$.
\end{definition}

This definition is justified by the definitions \ref{dkf} and \ref{norm}. When $Df(z)$ is invertible then, by analyticity, $\ga(f,z)$ is finite. We also have to consider the following number $K_\ze$ related to the sectional 
curvature at $\ze \in \MM_n$.

\begin{definition} \label{def-K} For any $\ze \in \MM_n$ 
$$ K_\ze = \sup \frac{d(\exp_z(u),\exp_z(v))}{\| u-v \|_z}$$
where the supremum is taken for all $z \in B_{\MM_n}(\ze , {\bf r}_\ze)$, and 
$u$, $v \in T_z\MM_n$ with $\| u \|_z$ and $\| v \|_z \leq {\bf r}_\ze),$ with
${\bf r}_\ze$ the radius of injectivity at $\ze$.
\end{definition}

\begin{remark} \label{rem1}
\begin{itemize}
\item $ K_\ze$ measures how fast the geodesics spread apart in $\MM_n$. 
When $u=0$ or more generally when $u$ and $v$ are on the same line through $0$,
$${d(\exp_z(u),\exp_z(v))}={\| u-v \|_z}.$$
Therefore, we always have
$$K_\ze \ge 1.$$
\item When $\MM_n$ has non-negative sectional curvature, the geodesics spread apart less than the rays (Do Carmo, \cite{Docarmo} Chap. V-2) so that 
$${d(\exp_z(u),\exp_z(v))} \le {\| u-v \|_z} $$
and consequently
$$K_\ze = 1.$$
\item Examples of manifolds with non-negative curvature are given by $\R^n$, $\S^n$ the unit sphere in $\R^{n+1}$, $\P^n(\R)$ the real projective space i.e the space of real vector lines in $\R^{n+1}$ (\cite{Docarmo}, Chap. 8, Prop. 4.4), $\P^n(\C)$ the complex projective space i.e the space of complex vector lines in $\C^{n+1}$ (\cite{Docarmo}, Chap. 8, Exerc. 11), a Lie group with a bi-invariant metric (\cite{Docarmo}, Chap. 4, Exerc. 1), $\O_n$ and $\S\O_n$ the orthogonal and special orthogonal groups (Lie groups) \ldots 
\end{itemize}
\end{remark}

\subsection{Main results for mappings.}

Our first main theorem relates the size of the quadratic attraction basin
of a zero $\ze$ of $f$ to the invariants $\ga(f,\ze)$ and $K_\ze$.

\begin{theorem} \label{th3} {\bf (R$-\ga-$theorem)} Let $f: \MM_n \rightarrow \R^n$ be analytic. 
Suppose that $f(\ze) = 0$ and $Df(\ze)$ is an isomorphism. Let
$$R(f,\ze) = \min \left ( {\bf r}_\ze , \frac{K_\ze + 2 - 
\sqrt{K_\ze^2 + 4K_\ze +2}}{2\ga(f,\ze)}\right ).$$
If $d(z, \ze) \leq R(f,\ze)$
then the Newton sequence $z_k = N_f^{(k)} (z)$ is defined for all $k \geq 0$, 
and 
$$ d(z_k, \ze) \leq \left ( \frac{1}{2} \right )^{2^k - 1}d(z, \ze).$$
\end{theorem}

\begin{remark}When $\MM_n = \R^n$ equipped with the usual metric structure, the
radius of injectivity ${\bf r}_\ze = \infty$ and $K_\ze = 1$. Thus, $R(f,\ze) = 
(3 - \sqrt{7})/
2\ga(f,\ze)$ as in Theorem \ref{th1}.

When $\MM_n$ has non-negative sectional curvature, according to Remark \ref{rem1} one has $K_\ze = 1$ and Theorem 
\ref{th3} becomes
\end{remark}

\begin{corollary} \label{corthe3} When $\MM_n$ has non-negative sectional curvature, let $f: \MM_n \rightarrow \R^n$ be analytic.
Suppose that $f(\ze) = 0$ and $Df(\ze)$ is an isomorphism. Let
$$R(f,\ze) = \min \left ( {\bf r}_\ze , \frac{3 - \sqrt{7}}{2\ga(f,\ze)}\right ).$$
If $d(z, \ze) \leq R(f,\ze)$
then the Newton sequence $z_k = N_f^{(k)} (z)$ is defined for all $k \geq 0$, 
and 
$$ d(z_k, \ze) \leq \left ( \frac{1}{2} \right )^{2^k - 1}d(z, \ze).$$
\end{corollary}

Theorem \ref{th3} has two interesting and immediate consequences: a lower estimate
for the distance from other zeros and a lower estimate
for the distance from the singular locus 
$$\Sigma_f = \{ z \in \MM_n \ : \ \det Df(z) = 0 \}.$$

\begin{corollary} \label{cor1} Suppose that $f(\ze) = 0$ and
$Df(\ze)$ is an isomorphism. Then, for any other zero $\ze' \not = \ze$ one has
$$d(\ze' , \ze) > R(f,\ze).$$
Moreover, for any $z \in \Sigma_f$ the same inequality holds:
$$d(z , \ze) > R(f,\ze).$$
\end{corollary}

Our second main theorem generalizes Theorem \ref{th2}. 
We give sufficient conditions for $z \in \MM_n$ to be the starting point
of a quadratically convergent Newton sequence. These conditions are
given in terms of $f$ at $z$, not in the behaviour of $f$ in a neighborhood of
$z$ as in Kantorovich theory. We first need three definitions.

\begin{definition} The function $\psi(u) = 1 - 4u + 2u^2$ is decreasing from $1$ to
$0$ when $0 \leq u \leq 1 - \sqrt{2}/2$. We denote by $\alpha_0 
= 0.130716944 \ldots $ the unique root of the equation $2u = \psi(u)^2$ in this 
interval. 
\end{definition}

\begin{definition} $\sigma$ is the sum of the following series:
$$\sigma = \sum_{k \geq 0} \left ( \frac{1}{2} \right )^{2^k - 1} = 1.632843018 \ldots $$
\end{definition}

\begin{definition} 
$$s_0 = \frac{1}{\sigma + \frac{(1-\sigma\alpha_0)^2}{\psi(\sigma\alpha_0)}\left(1+\frac{\sigma}{1-\sigma\alpha_0}\right)} = 0.103621842 \ldots 
$$
\end{definition}

\begin{definition} We let $\be (f,z) = \| Df(z)^{-1}f(z) \|_z$ and 
$\al (f,z) = \be (f,z) \ga (f,z)$. We give to $\be (f,z)$ and $\al (f,z)$
the value $\infty$ when $Df(z)$ is singular.
\end{definition}

\begin{theorem} \label{th4} {\bf (R$-\alpha-$Theorem)} Let $f : \MM_n \rightarrow \R^n$ be analytic. 
Let $\ {z \in \MM_n}$ be such that
$$\be (f,z) \leq s_0 {\bf r}_z \ \ \mbox{and} \ \ \al (f,z) < \al_0.$$ 
Then the Newton sequence $z_0 = z$, $z_{k+1} = N_f(z_k)$
is defined for all integers $k \geq 0$ and converges to a zero $\ze$ of $f$. Moreover,
$$d(z_{k+1}, z_k) \leq \left ( \frac{1}{2} \right )^{2^k - 1}\be (f,z)$$
and 
$$d(\ze , z) \leq \sigma \be (f,z).$$
\end{theorem}

\begin{remark}When $\MM_n = \R^n$ is equipped with the usual metric structure, the
radius of injectivity ${\bf r}_\ze = \infty$ and the first condition in 
Theorem \ref{th4} is automatically satisfied. In this context Theorems 
\ref{th2} and \ref{th4} coincide.
\end{remark}

\subsection{Main results for vector fields.}

The case of vector fields is treated similarly. As in Theorem \ref{th3} we have:

\begin{theorem} \label{thV3} {\bf (R$-\ga-$Theorem)} Let $X: \MM_n \rightarrow T\MM_n$ be an analytic vector field. 
Suppose that $X(\ze) = 0$ and $DX(\ze)$ is an isomorphism. Let
$$R(X,\ze) = \min \left ( {\bf r}_\ze , \frac{K_\ze + 2 - \sqrt{K_\ze^2 + 4K_\ze +2}}{2\ga(X,\ze)}\right ).$$
If $d(z, \ze) \leq R(X,\ze)$
then the Newton sequence $z_k = N_X^{(k)} (z)$ is defined for all $k \geq 0$, 
and 
$$ d(z_k, \ze) \leq \left ( \frac{1}{2} \right )^{2^k - 1}d(z, \ze).$$
\end{theorem}

Like for mappings, Theorem \ref{thV3} gives estimates
for the distance from other zeros and a lower estimate
for the distance from the singular locus 
$$\Sigma_X = \{ z \in \MM_n \ : \ \det DX(z) = 0 \}.$$

\begin{corollary} \label{corV1} Suppose that $X(\ze) = 0$ and
$DX(\ze)$ is an isomorphism. Then, for any other zero $\ze' \not = \ze$ one has
$$d(\ze' , \ze) > R(X,\ze).$$
Moreover, for any $z \in \Sigma_X$ the same inequality hold:
$$d(z , \ze) > R(X,\ze).$$
\end{corollary}

The invariants $\be$ and $\al$ are defined similarly:

\begin{definition} We let 
$$\be (X,z) = \| DX(z)^{-1}X(z) \|_z$$ 
and 
$$\al (X,z) = \be (X,z) \ga (X,z).$$ 
We give to $\be (X,z)$ and $\al (X,z)$
the value $\infty$ when $DX(z)$ is singular.
\end{definition}

\begin{theorem} \label{thV4} {\bf (R$-\a-$Theorem)} Let $X : \MM_n \rightarrow T\MM_n$ be an analytic vector field. Let $z \in \MM_n$ be such that
$$\be (X,z) \leq s_0 {\bf r}_z \ \ \mbox{and} \ \ \al (X,z) < \al_0.$$ 
Then the Newton sequence $z_0 = z$, $z_{k+1} = N_X(z_k)$
is defined for all integers $k \geq 0$ and converges to a zero $\ze$ of $X$. Moreover,
$$d(z_{k+1}, z_k) \leq \left ( \frac{1}{2} \right )^{2^k - 1}\be (X,z)$$
and 
$$d(\ze , z) \leq \sigma \be (X,z).$$
\end{theorem}

\subsection{Previous work.}

There is quite a bit of previous work on such questions. The first to consider
Newton's method on a manifold is Rayleigh 1899 \cite{Rayleigh} who
defined what we call today ``Rayleigh Quotient Iteration'' which is
in fact a Newton iteration for a vector field on the sphere. Then, Shub 1986 
\cite{Shub-Venezuela} defined Newton's method for the problem of
finding the zeros of a vector field on a manifold and used
retractions to send a neighborhood of the origin in the tangent space
onto the manifold itself. In our paper we do not use general
retractions but exponential maps. 
Independently of~\cite{Shub-Venezuela}, Smith 1994
\cite{Smith} developed an intrinsic Newton's method and a conjugate
gradient algorithm on a manifold using the exponential map. 
Also independently,
Udriste 1994 \cite{Udriste} studied 
Newton's method to find the zeros of a gradient vector field defined 
on a Riemannian manifold; Owren and Welfert 1996 \cite{Owren} defined Newton's
iteration for solving the equation $F(x)=0$ where $F$ is a map
from a Lie group to its corresponding Lie algebra; 
%Smith 1994 \cite{Smith} and 
Edelman-Arias-Smith 1998 \cite{Edelman-Arias-Smith}
developed Newton's and conjugate gradient algorithms on the
Grassmann and Stiefel manifolds. These authors define
Newton's method via the exponential map as we do here. 
Shub 1993 \cite{Shub93}, Shub and Smale 1993-1996 \cite{Bez1}, \cite{Bez2}, 
\cite{Bez3}, \cite{Bez4}, \cite{Bez5}, see also, Blum-Cucker-Shub-Smale 1998 
\cite{BCSS}, Malajovich 1994 \cite{Malajovich},
Dedieu and Shub 2000 \cite{Dedieu-Shub-2} introduce and study Newton's
method on projective spaces and their products. Another important paper about
this subject is Adler-Dedieu-Margulies-Martens-Shub 2001 \cite{ADMMS} where
qualitative aspects of Newton's method on Riemannian manifolds are investigated 
for both mappings and vector fields. This paper contains a nice application 
to a geometric model for the human spine represented as a $18-$tuple of 
$3 \times 3$ orthogonal matrices. Recently Ferreira-Svaiter \cite{Benar} gave 
a Kantorovich-like theorem for Newton's method for vector fields defined on 
Riemannian manifolds. 

\section{Parallel transport and Taylor's formula.}

In the proof sections of this paper, we frequently use parallel transport:

\begin{definition} {\bf (Parallel transport.)} Let $z_0$ and $z \in \MM_n$ with $z$ in the ball about $z_0$
with radius ${\bf r}_{z_0}$ the radius of injectivity. Then, there exists a unique geodesic curve $c(t)$ in 
this ball such that $c(0) = z_0$ and $c(T) = z$ for a certain $T$. In this context we denote by
$$P_{z_0,z} : T_{z_0}\MM_n \rightarrow T_z\MM_n$$
the parallel transport along this geodesic. It is an isometry which preserves the orientation when 
$\MM_n$ is oriented.
\end{definition} 

We now extend this concept to other objects

\begin{definition} {\bf (Parallel transport: extension.)} 
\begin{itemize}
\item For a covector $\omega_{z_0} \in T_{z_0}^*\MM_n$ by
$$(P_{z_0,z}\omega_{z_0})(Y_z) = \omega_{z_0}(P_{z_0,z}^{-1}(Y_z))$$
for any $Y_z \in T_z\MM_n$. 
\item For a tensor field $T \in \T_q^p(\MM_n)$ we denote by $T_{z_0}$ its value at ${z_0}$ that is 
$$T_{z_0}((\omega^1)_{z_0} \ldots (\omega^p)_{z_0},(Y_1)_{z_0} \ldots (Y_q)_{z_0})=
T(\omega^1 \ldots \omega^p,Y_1 \ldots Y_q)(z_0)$$
for any $1-$forms $\omega^i$ and vector fields $Y_j$.
\item Parallel transport for $T_{z_0}$ is defined by
$$P_{z_0,z}T_{z_0}((\omega^1)_z \ldots (\omega^p)_z,(Y_1)_z \ldots (Y_q)_z) = $$
$$T_{z_0}(P_{z_0,z}^{-1}(\omega^1)_z \ldots P_{z_0,z}^{-1}(\omega^p)_z,
P_{z_0,z}^{-1}(Y_1)_z \ldots P_{z_0,z}^{-1}(Y_q)_z)$$
for any covectors $(\omega^i)_z \in T_{z}^*\MM_n$ and vectors $(Y_j)_z \in T_{z}\MM_n$.
\end{itemize}
\end{definition}

The covariant derivative of a tensor field $T \in \T_q^p(\MM_n, \R^m)$ at a point 
may be described in terms of parallel transport: 
let $z_0 \in \MM_n$ and $u \in T_{z_0}\MM_n$ be given. With the geodesic curve $c(t) = exp_{z_0}(tu)$ we have
$$DT(z_0)u = \lim_{t \rightarrow 0} \frac{1}{t} \left ( P_{z_0,c(t)}^{-1}T_{c(t)} - T_{z_0} \right ).$$

We now give Taylor's formula. A reference is Dieudonn\'e \cite{Dieud}, Chap. XVIII-6, where the case of functions is considered. Tensors are treated similarly. We have

\begin{theorem} \label{th5} {\bf (Taylor formula)} For any tensor field $T \in \T_q^p(\MM_n, \R^m)$,
$z_0$, $z \in \MM_n$ with $z$ in a certain neighborhood about $z_0$, and $u \in T_{z_0}\MM_n$ such that
$z = \exp_{z_0}(u)$ we have
$$T(z) = \left ( \sum_{k=0}^\infty \frac{1}{k!} D^kT(z_0)u^k
\right )P_{z,z_0}.$$
\end{theorem}

Taking the $l-$th covariant derivative in \ref{th5} gives the following:

\begin{corollary} \label{cth5} With the same hypothesis, for any $l \ge 0$, we have
$$D^lT(z) = \left ( \sum_{k=0}^\infty \frac{1}{k!} D^{k+l}T(z_0)u^k 
\right )P_{{z,z_0}}.$$
\end{corollary}

The neighborhood of $z_0$ in Theorem~\ref{th5} and Corollary~\ref{cth5} is given by the radius of injectivity at $z_0$ and by the disks 
of convergence of the Taylor series of the coordinates of the tensor field $T$ in a local chart about $z_0$. 
In the following we relate it to $\ga(f,z).$
Let $f : \MM_n \rightarrow \R^n$ be an analytic map (resp. $X : \MM_n \rightarrow T\MM_n$ an analytic vector field). 
As an immediate consequence of the definition of $\ga(f,z)$ (resp. $\ga (X,z))$ we have:

\begin{proposition} The Taylor series at $z \in \MM_n$ for $f$ and $D^kf$ (resp. $X$ and $D^kX$) converge in the ball about $z$ with radius $1/\ga(f,z)$ (resp. $1/\ga(X,z)$). Theorem \ref{th5} is valid for any $z$ with 
$$d(z,z_0) < \min ({\bf r}_{z_0}, 1/\ga(f,z)) \ \ \mbox{(resp.} \ \ \min ({\bf r}_{z_0}, 1/\ga(X,z))).$$
\end{proposition}

\noindent {\bf Proof.} Taking a local chart it suffices to prove this theorem in the context
of a map $f : \R^n \rightarrow \R^n$. Then, we use \cite{BCSS} Chap. 8, Prop. 6. \qed

\section{Proof of the R$-\ga-$theorem.}

This proof is quite long and split in a series of lemmas. We frequently use the notations
$\| A \|_{E,F}$ for the operator norm of the linear map $A : E \rightarrow F$ and
$\| A \|_{E}$ when $E=F$.

\begin{lemma} \label{lem1} Let $x$, $y\in \MM_n$ with $d(x,y) < {\bf r}_x$. We suppose
that $Df(x)$ is non-singular and that 
$$Df(x)^{-1}Df(y) = P_{y,x} + BP_{y,x}$$
with $\| B \|_{T_x\MM_n} \leq r$ for a certain $r<1$. Then, $Df(y)$ is non-singular and
$$\| Df(y)^{-1}Df(x) \|_{T_x\MM_n,T_y\MM_n} \leq \frac{1}{1-r}.$$
\end{lemma}

\noindent {\bf Proof.} $Df(x)^{-1}Df(y) = (id_{T_x\MM_n} + B)P_{y,x}$. Since 
$\| B \|_{T_x\MM_n} \leq r < 1$ the operator $id_{T_x\MM_n} + B$ is non-singular and
its inverse satisfies $ \| id_{T_x\MM_n} + B \|_{T_x\MM_n} \| \leq 1/(1 - r).$ Then, we notice
that parallel transport $P_{y,x}$ is an isometry. \qed

\begin{lemma}\label{lem2}Let $x$, $y\in \MM_n$ with $d(x,y) < {\bf r}_x$. We suppose that
$Df(x)$ is non-singular and that 
$$\nu = d(x,y)\ga(f,x) < 1 - \frac{\sqrt{2}}{2}.$$
Then, $Df(y)$ is non-singular and
$$\| Df(y)^{-1}Df(x) \|_{T_x\MM_n,T_y\MM_n} \leq \frac{(1 - \nu)^2}{\psi(\nu)}.$$
\end{lemma}

\noindent {\bf Proof.} Let $u = \exp_x^{-1}(y)$. By Corollary \ref{cth5}
with $l=1$ and $T=f$ we get
$$Df(y) = \left ( \sum_{k=0}^\infty \frac{1}{k!} D^{k+1}f(x)u^k \right ) P_{y,x}.$$
so that
$$Df(x)^{-1}Df(y) = P_{y,x} + \left ( \sum_{k=2}^\infty \frac{1}{(k-1)!} 
Df(x)^{-1}D^{k}f(x)u^{k-1} \right ) P_{y,x} = $$
$$P_{y,x} + BP_{y,x}.$$
Let us now give a bound for $\| B \|_{T_x\MM_n}$:
$$\| B \|_{T_x\MM_n} = \left \| \sum_{k=2}^\infty \frac{1}{(k-1)!} Df(x)^{-1}D^{k}f(x)
u^{k-1} \right \|_{T_x\MM_n} \leq $$
$$\sum_{k=2}^\infty k \frac{1}{k!}\| Df(x)^{-1}D^{k}f(x) \|_{T_x\MM_n} 
\| u \|_x^{k-1} \leq \sum_{k=2}^\infty k \ga(f,x)^{k-1}\| u \|_x^{k-1} = $$
$$\sum_{k=2}^\infty k \ga(f,x)^{k-1} d(x,y)^{k-1} = \sum_{k=2}^\infty k \nu^{k-1} = 
\frac{1}{(1 - \nu)^2} - 1.$$
This last quantity is $<1$ because $\nu < 1 - \frac{\sqrt{2}}{2}$. The conclusion is obtained 
from Lemma \ref{lem1}. \qed

\begin{lemma}\label{lem3}Let $z$, $\ze \in \MM_n$ with $d(z,\ze) < {\bf r}_\ze $. 
We suppose that
$f(\ze)=0$, $Df(\ze)$ is non-singular and 
$$\nu = d(z,\ze)\ga(f,\ze) < 1 - \frac{\sqrt{2}}{2}.$$
Then, 
$$\| Df(\ze)^{-1}(Df(z)\exp_z^{-1}(\ze) + f(z)) \|_\ze \leq 
\frac{ \nu d(z,\ze)}{(1-\nu)^2}.$$
\end{lemma}

\begin{remark} Let $u \in T_\ze \MM_n$ be such that $\exp_\ze (u) = z$.
Let $v = P_{\ze, z}u \in T_z \MM_n$ be the parallel transport of $u$ along the geodesic
between $\ze$ and $z$. Then, $\exp_z(-v) = \ze$ so that the expression the vector
$\exp_z^{-1}(\ze)$ is equal to $-u$.
\end{remark}

\noindent {\bf Proof.} Let $u \in T_\ze \MM_n$ be such that $\exp_\ze (u) = z$. From Taylor formula
we get:
$$f(z) = f(\ze) + Df(\ze)u + \sum_{k \geq 2} \frac{1}{k!}D^kf(\ze)u^k$$
and
$$Df(z) = \left ( Df(\ze) + \sum_{k \geq 2} \frac{k}{k!}D^kf(\ze)u^{k-1} \right )P_{z,\ze}.$$
Notice $f(\ze) = 0$ and $P_{z,\ze}\exp_z^{-1}(\ze) = -u$ thus,
$$Df(\ze)^{-1}(Df(z)\exp_z^{-1}(\ze) + f(z)) = - \sum_{k \geq 2} \frac{k-1}{k!}
Df(\ze)^{-1}D^kf(\ze)u^k$$
and
$$\| Df(\ze)^{-1}(Df(z)\exp_z^{-1}(\ze) + f(z)) \|_\ze \leq \sum_{k \geq 2} (k-1)
\ga(f,\ze)^{k-1} \| u \|_\ze^k =$$
$$\sum_{k \geq 2} (k-1) \ga(f,\ze)^{k-1} d(z,\ze)^k = \frac{ \nu d(z,\ze)}{(1-\nu)^2}$$
and we are done. \qed

\begin{lemma}\label{lem4}Let $z$, $\ze \in \MM_n$ with $d(z,\ze) < {\bf r}_\ze $. 
We suppose that
$f(\ze)=0$, $Df(\ze)$ is non-singular and 
$$\nu = d(z,\ze)\ga(f,\ze) < 1 - \frac{\sqrt{2}}{2}.$$
Then, 
$$\| \exp_z^{-1}(N_f(z)) - \exp_z^{-1}(\ze) \|_z \leq \frac{ \nu d(z,\ze)}{\psi(\nu)}.$$
\end{lemma}

\noindent {\bf Proof.} $$\| \exp_z^{-1}(N_f(z)) - \exp_z^{-1}(\ze)\|_z = 
\| Df(z)^{-1}(Df(z)\exp_z^{-1}(\ze) + f(z))\|_z = $$
$$\| Df(z)^{-1}Df(\ze)Df(\ze)^{-1}(Df(z)\exp_z^{-1}(\ze) + f(z))\|_z \leq $$
$$\| Df(z)^{-1}Df(\ze)\|_{T_{\ze}\MM_n,T_z \MM_n}
\| Df(\ze)^{-1}(Df(z)\exp_z^{-1}(\ze) + f(z))\|_\ze \leq \frac{(1 - \nu)^2}{\psi(\nu)}
\frac{ \nu d(z,\ze)}{(1-\nu)^2}$$
by Lemma \ref{lem2} and Lemma \ref{lem3}. This achieves the proof. \qed

Let us recall the definition of the geometric constant 
$$ K_\ze = \sup \frac{d(\exp_z(u),\exp_z(v))}{\| u-v \|_z}$$
where the supremum is taken for all $z \in B_{\MM_n}(\ze , {\bf r}_\ze)$, and 
$u$, $v \in T_z\MM_n$ with $\| u \|_z$ and $\| v \|_z \leq {\bf r}_\ze).$

\begin{lemma}\label{lem5} The following inequalities hold: $ K_\ze \geq 1$ and 
$$\frac{K_\ze+2 - \sqrt{K_\ze^2 + 4K_\ze + 2}}{2} < 1 - \frac{\sqrt{2}}{2}.$$ Moreover, if 
$$\nu \leq \frac{K_\ze+2 - \sqrt{K_\ze^2 + 4K_\ze + 2}}{2}$$ then 
$$\frac{K_\ze\nu }{ \psi(\nu)} \leq \frac{1}{2}.$$
\end{lemma}

\noindent {\bf Proof.} The constant $K_\ze$ is necessarily $\geq 1$ because 
$$d(\exp_\ze(0), \exp_\ze(v)) = \| 0-v \|_\ze$$ for any $v$ in the ball
of injectivity for $\ze$. The second inequality comes from
$$\frac{K_\ze+2 - \sqrt{K_\ze^2 + 4K_\ze + 2}}{2} = \frac{1}{K_\ze+2 + \sqrt{K_\ze^2 + 4K_\ze + 2}} \leq
\frac{1}{3 + \sqrt{7}} < 1 - \frac{\sqrt{2}}{2}.$$
The third inequality uses the fact ${K_\ze\nu }/{ \psi(\nu)}$ is increasing on 
the interval $[0,1 - \frac{\sqrt{2}}{2}[.$ \qed

\vskip 3mm
\noindent {\bf Proof of Theorem \ref{th3}.} We are going to prove that
$$d(z_k, \ze) \leq \left ( \frac{K_\ze \nu}{\psi(\nu)} \right )^{2^k - 1} d(z, \ze)$$
for any $k \geq 0$ with $\nu = d(z,\ze)\ga(f,\ze)$. The conclusion is then an easy
consequence of the hypothesis and of Lemma \ref{lem5}. We proceed by induction:
the case $k=0$ is evident. Then,
$$d(z_{k+1}, \ze) = d(N_f(z_{k}), \ze) \leq K_\ze \| \exp_{z_k}^{-1}(N_f(z_k)) - 
\exp_{z_k}^{-1}(\ze) \|_{z_k}.$$
>From Lemma \ref{lem4} we get 
$$d(z_{k+1}, \ze) \leq K_\ze \frac{ \nu_k d(z_k,\ze)}{\psi(\nu_k)} $$
with $\nu_k = d(z_k,\ze)\ga(f,\ze)$. By the induction hypothesis
$$d(z_{k+1}, \ze) \leq \frac{ K_\ze \ga(f,\ze)}{\psi(\nu)}
\left ( \left ( \frac{K_\ze \nu}{\psi(\nu)} \right )^{2^k - 1} d(z, \ze)\right )^2 = 
\left ( \frac{K_\ze \nu}{\psi(\nu)} \right )^{2^{k+1} - 1} d(z, \ze).$$
\qed

\section{Proof of the R$-\al-$theorem.}

Let us first recall two definitions: 
$\be(f,z) = \| Df(z)^{-1}f(z) \|_z$ and $\a = \be \ga$ . For the proof of Theorem \ref{th4}
with need some more lemmas. 

\begin{lemma}\label{lem6} For $\vert r \vert < 1$ and any integer $k \geq 0$
$$\sum_{l=0}^\infty \frac{(k+l)!}{k!\ l!}\ r^l = \frac{1}{(1-r)^{k+1}}.$$
\end{lemma}

\begin{lemma}\label{lem7}Let $z$, $z_1 \in \MM_n$ with $d(z,z_1) < {\bf r}_z $. 
We suppose that $Df(z)$ is nonsingular and
$$\nu = d(z,z_1)\ga(f,z) < 1 - \frac{\sqrt{2}}{2}.$$
Then, for any integer $k \geq 2$
\begin{itemize}
\item
$\ds\frac{\| Df(z_1)^{-1}D^kf(z_1) \|_{z_1}}{k!} \leq \ds\frac{1}{\psi(\nu)}
\left ( \ds\frac{\ga(f,z)}{1-\nu} \right )^{k-1},$
\item
$\| Df(z)^{-1}f(z_1) \|_{z} \leq \be(f,z) + \ds\frac{d(z,z_1)}{1-\nu}.$
\end{itemize}
\end{lemma}

\noindent {\bf Proof.} Let $u \in T_\ze \MM_n$ be such that $\exp_\ze (u) = z$. 
>From Taylor formula (Theorem \ref{cth5}) we get:
$$D^kf(z_1) = \left ( \sum_{l \geq 0} \frac{1}{l!}D^{k+l}f(z)u^l
\right )P_{z,z_1}$$
so that
$$\frac{\| Df(z_1)^{-1}D^kf(z_1) \|_{z_1}}{k!} \leq $$
$$\| Df(z_1)^{-1}Df(z) \|_{T_{z_1}\MM_n, T_{z}\MM_n} \left \| \left ( \sum_{l \geq 0} 
\frac{1}{k! \ l!}Df(z)^{-1}D^{k+l}f(z)u^l \right )P_{z, z_1} \right \|_{z_1} \leq $$
$$\frac{(1-\nu)^2}{\psi(\nu)} \sum_{l \geq 0} \frac{(k+l)!}{k! \ l!}\ga(f,z)^{k+l-1}
\| u \|_z^l = \frac{(1-\nu)^2}{\psi(\nu)} \ga(f,z)^{k-1} \frac{1}{(1-\nu)^{k+1}}$$
using Lemma \ref{lem2} and Lemma~\ref{lem6}, the definition of $\ga$, 
the fact that $P_{z,z_1}$ is an isometry
and $\| u \|_z = d(z,z_1)$. This proves the first inequality. Let us now prove 
the second one.
$$\| Df(z)^{-1}f(z_1) \|_{z_1} = \left \| \sum_{k \geq 0} \frac{1}{k!}Df(z)^{-1}
D^{k}f(z)u^k P_{z,z_1} \right \|_{z_1} \leq $$
$$\| Df(z)^{-1}f(z)\|_{z} + \| u \|_z \sum_{k \geq 1} \ga(f,z)^{k-1} 
\| u \|_z^{k-1} = \be(f,z) + d(z,z_1)\frac{1}{1 - \nu}$$
and we are done. \qed

\begin{lemma}\label{lem8} Let $z$, $z_1 \in \MM_n$ with $d(z,z_1) < {\bf r}_z $. 
We suppose that $Df(z)$ is nonsingular and
$$\nu = d(z,z_1)\ga(f,z) < 1 - \frac{\sqrt{2}}{2}$$
then
\begin{itemize}
\item $\be(f,z_1) \leq \displaystyle\frac{(1-\nu)^2}{\psi(\nu)}\left(\be(f,z) + \displaystyle\frac{d(z,z_1)}{1 - \nu}\right),$
\item$\ga(f,z_1) \leq \displaystyle\frac{\ga(f,z)}{(1-\nu)\psi(\nu)}.$
\end{itemize}
\end{lemma}

\noindent {\bf Proof.} 
The first estimate is a consequence of Lemma \ref{lem2}, Lemma \ref{lem7} and the following
$$\be(f,z_1) = \left \| Df(z_1)^{-1} f(z_1) \right \|_{z_1} \le \left \| Df(z_1)^{-1} Df(z) \right \|_{z, z_1} \left \| Df(z)^{-1} f(z_1) \right \|_{z}.$$
The second inequality is an easy consequence of Lemma \ref{lem7}
$$\ga(f,z_1) = \sup_{k \geq 2} \left \| Df(z_1)^{-1} \frac{D^kf(z_1)}{k!}
\right \|^{1/k-1} \leq \frac{\ga(f,z)}{1-\nu} \sup_{k \geq 2}\frac{1}{\psi(\nu)^{1/k-1}} 
= \frac{\ga(f,z)}{(1-\nu)\psi(\nu)}$$
because $\nu < 1 - \sqrt{2}/2$ implies $\psi(\nu) < 1$ and the supremum is achieved 
for $k=2$. 
\qed

\begin{lemma}\label{lemn1} Let $\MM_n$ be a complete Riemannian manifold. Then, for any $x$, $y \in \MM_n$ we have
$${\bf r}_x - d(x,y) \le {\bf r}_y.$$
\end{lemma}

\noindent {\bf Proof.} To prove this inequality we show that $\exp_y$ is injective in the ball about $0$ with radius ${\bf r}_x - d(x,y)$ in $T_y\MM_n$. Let $u \in T_x\MM_n$ be such that $y = \exp_x u$ and $\| u \|_x = d(x,y).$ Let $v$, $w \in T_y\MM_n$ be such that $\exp_y (v) = \exp_y (w)$ and $\| v \|_y = \| w \|_y < {\bf r}_x - d(x,y).$ Let $P$ denote the parallel transport from $T_y\MM_n$ to $T_x\MM_n$. We have
$$\exp_y (v) = \exp_x (u + Pv) \ \ \mbox{and} \ \ \exp_y (w) = \exp_x (u + Pw).$$
Moreover 
$$\| u + Pv \|_x \le \| u \|_x + \| Pv \|_x = d(x,y) + \| v \|_y < d(x,y) + {\bf r}_x - d(x,y) = {\bf r}_x$$
and a similar inequality holds with $w$. Since $\exp_x$ is injective in this ball we get $u + Pv = u + Pw$ so that $v=w$ and we are done. \qed

\begin{lemma}\label{lemn2} Let $x \in \MM_n$ be such that
$$\be (f,x) \leq s_0 {\bf r}_x \ \ \mbox{and} \ \ \al (f,x) < \al_0.$$ 
Then, for any $y \in \MM_n$ such that $d(x,y) \le \sigma \be (f,x)$ we have
$$\be(f,y) \le {\bf r}_y.$$
\end{lemma}

\noindent {\bf Proof.} Let $s$ be a positive real number and let us suppose that $\be (f,x) \leq s {\bf r}_x.$ Let $y$ be such that $d(x,y) \le \sigma \be (f,x)$. We have, by Lemma \ref{lem8} 
$$\be(f,y) \leq \frac{(1-\nu)^2}{\psi(\nu)}\left(\be(f,x) + \displaystyle\frac{d(x,y)}{1 - \nu}\right) \le$$
$$\frac{(1-\nu)^2}{\psi(\nu)}\left(1 + \frac{\sigma}{1-\nu}\right)\be(f,x) \le 
\frac{(1-\nu)^2}{\psi(\nu)}\left(1 + \frac{\sigma}{1-\nu}\right)s{\bf r}_x.$$
Moreover, by Lemma \ref{lemn1},
$${\bf r}_x \le {\bf r}_y + d(x,y) \le {\bf r}_y + \sigma \beta(f,x) \le {\bf r}_y + \sigma s {\bf r}_x$$
so that 
$${\bf r}_x \le \frac{1}{1 - \sigma s}{\bf r}_y$$
as soon as $\sigma s < 1$. Thus
$$\be(f,y) \leq \frac{(1-\nu)^2}{\psi(\nu)}\left(1 + \frac{\sigma}{1-\nu}\right)\frac{s}{1 - \sigma s}{\bf r}_y$$
so that $\be(f,y) \leq {\bf r}_y$ if
$$\frac{(1-\nu)^2}{\psi(\nu)}\left(1 + \frac{\sigma}{1-\nu}\right)\frac{s}{1 - \sigma s} < 1 \ \ \mbox{and} \ \ 
\sigma s < 1.$$
These conditions are satisfied when 
$$s \le \frac{1}{\sigma + \frac{(1-\nu)^2}{\psi(\nu)}\left(1 + \frac{\sigma}{1-\nu}\right)}.$$
We also notice that
$$\nu = d(x,y) \ga(f,x) \le \sigma \beta(f,x) \gamma (f,x) = \sigma \alpha(f,x) \le \sigma \alpha_0.$$
Since the function
$$\nu \rightarrow \sigma + \frac{(1-\nu)^2}{\psi(\nu)}\left(1 + \frac{\sigma}{1-\nu}\right)$$
is increasing we obtain the following sufficient condition
$$s \le \frac{1}{\sigma + \frac{(1-\sigma \alpha_0)^2}{\psi(\sigma \alpha_0)}\left(1 + \frac{\sigma}{1-\sigma \alpha_0}\right)} = 0.103621842 \ldots $$
\qed

\begin{lemma}\label{lem9} Let $z \in \MM_n$ and $z_1 = N_f(z)$. 
We suppose that
$$\nu = d(z,z_1)\ga(f,z) < 1 - \frac{\sqrt{2}}{2}$$
then
\begin{itemize}
\item
$\be(f,z_1) \leq \frac{1-\nu}{\psi(\nu)}\be(f,z)^2 \ga(f,z),$
\item
$\a(f,z_1) \leq \frac{\a(f,z)^2}{\psi(\nu)^2}.$
\end{itemize}
\end{lemma}

\noindent {\bf Proof.} From Lemma \ref{lem2} we get the following
$$\be(f,z_1) = \| Df(z_1)^{-1}f(z_1) \|_{z_1} \leq \| Df(z_1)^{-1}Df(z) 
\|_{T_z\MM_n,T_{z_1}\MM_n} \| Df(z)^{-1}f(z_1) \|_z \leq $$
$$\frac{(1-\nu)^2}{\psi(\nu)}\| Df(z)^{-1}f(z_1) \|_z.$$
Let $u \in T_z\MM_n$ be such that $\exp_z(u) = z_1$. From Taylor formula
$$f(z_1) = f(z) + Df(z)u + \sum_{k \geq 2}\frac{1}{k!}D^kf(z)u^k.$$
Since $z_1 = N_f(z)$ we have $f(z) + Df(z)u=0$ so that
$$\| Df(z)^{-1}f(z_1) \|_z \leq \sum_{k \geq 2}\frac{1}{k!}\| Df(z)^{-1}D^kf(z) \|_z
\| u \|_z^k \leq $$
$$\sum_{k \geq 2}\ga(f,z)^{k-1} d(z,z_1)^k = \frac{\ga(f,z) d(z,z_1)^2}{1-\nu} = 
\frac{\ga(f,z) \be(f,z)^2}{1-\nu}.$$
This proves the first inequality. For the second we multiply together
$\be(f,z_1) \leq \frac{1-\nu}{\psi(\nu)}\be(f,z)^2 \ga(f,z),$ and $\ga(f,z_1) \leq \frac{\ga(f,z)}{(1-\nu)\psi(\nu)}$ obtained in Lemma \ref{lem8}. \qed

\vskip 3mm
\noindent {\bf Proof of Theorem \ref{th4}.} Let us first introduce some more notations: 
$z_k$ is the Newton sequence starting at $z_0 = z$, $\be_k = \be(f,z_k) = d(z_k,
z_{k+1})$, $\ga_k = \ga(f,z_k)$, $\a_k = \a(f,z_k) = \ga_k d(z_k, z_{k+1}) 
$ and ${\bf r}_k$ the radius of injectivity at $z_k$. We shall prove, by induction, 
the following:
\begin{itemize} \item $1_k:$ $\a_k \leq \left ( \frac{1}{2} \right )^{2^k - 1} \a_0,$
\item $2_k:$ $\be_k \leq \left ( \frac{1}{2} \right )^{2^k - 1} \be_0,$
\item $3_k:$ $\be_k \leq {\bf r}_k.$
\end{itemize} 

This will prove Theorem \ref{th4}. These inequalities are clearly satisfied when $k=0$.
To prove $1_{k+1}$ we use $1_k$ and Lemma \ref{lem9}:
$$\a_{k+1} \leq \frac{\a_k^2}{\psi(\a_k)^2} \leq \frac{1}{\psi(\a_k)^2}
\left ( \left ( \frac{1}{2} \right )^{2^k - 1} \a_0 \right )^2 \leq 
\frac{\a_0}{\psi(\a_k)^2}
\left ( \frac{1}{2} \right )^{2^{k+1} - 2} \a_0 \leq $$
$$\frac{\a_0}{\psi(\a_0)^2}
\left ( \frac{1}{2} \right )^{2^{k+1} - 2} \a_0 \leq 
\frac{1}{2}
\left ( \frac{1}{2} \right )^{2^{k+1} - 2} \a_0 = 
\left ( \frac{1}{2} \right )^{2^{k+1} - 1} \a_0.$$
To prove $2_{k+1}$ we use a similar argument: by Lemma \ref{lem9}
$$\b_{k+1} \leq \frac{1-\a_k}{\psi(\a_k)} \ga_k \b_k^2 = 
\frac{1-\a_k}{\psi(\a_k)} \a_k \b_k \leq \frac{1-\a_0}{\psi(\a_0)}
\left ( \frac{1}{2} \right )^{2^k - 1} \a_0
\left ( \frac{1}{2} \right )^{2^k - 1} \be_0.$$
Since $\frac{\a_0(1-\a_0)}{\psi(\a_0)} \leq \frac{1}{2}$ we obtain
$$\b_{k+1} \leq \frac{1}{2} \left ( \left ( \frac{1}{2} \right )^{2^k - 1}
\right )^2 \be_0 = \left ( \frac{1}{2} \right )^{2^{k+1} - 1} \be_0$$
and we are done. To prove $3_{k+1}$ we use the hypothesis, Lemma \ref{lemn2} and the following 
estimate:
$$d(z_{k+1},z_0) \leq \sum_{i=0}^{k}d(z_{i+1},z_i) \leq \sum_{i=0}^{k}\b_i
\leq \sum_{i=0}^{k}\left ( \frac{1}{2} \right )^{2^i - 1} \be_0 \leq \sigma \be_0.$$
\qed

%\section{Proofs for the Case of Vector Fields}

%They are similar ``mutatis mutandis'' to the proofs of Theorems \ref{th3} and \ref{th4}
%and based on Taylor Formula which is still valid in this new context.

\noindent {\bf Proof of Theorems \ref{thV3} and \ref{thV4}.} 
The proofs of Theorems \ref{thV3} and \ref{thV4} are formally identical
to the proofs of Theorems \ref{th3} and \ref{th4},
respectively. The only difference is that $X$ is an analytic vector field,
and vector fields are 1-contravariant 0-covariant tensor fields.
Therefore, its $k$-th derivative is a 1-contravariant $k$-covariant
tensor field, instead of a $k$-covariant tensorial vector field.

\section{Examples}\label{examples}

\noindent {\bf First example: the unit sphere. } $\S^{n}$ denotes the unit sphere in $\R^{n+1}$, the tangent space $T_x\S^{n}$ is the hyperplane in $\R^{n+1}$ orthogonal to $x$, the Riemannian structure is given by the Euclidean structure of $\R^{n+1}$ and the Riemannian distance in $\S^{n}$ is the arc length taken along great circles: $$d(x,y) = \arccos \langle x,y \rangle.$$
The exponential map at $x \in \S^{n}$ is given by
$$\exp_x (u) = x \cos \| u \| + u \frac{\sin \| u \|}{\| u \|}$$
for any $u \in T_x\S^{n}.$ The radius of injectivity is equal to ${\bf r}_x = \pi$ and the constant appearing in Definition \ref{def-K} is $K_x = 1$ because $\S^{n}$ has positive sectional curvature. 
Newton's method is given by 
\begin{eqnarray*}
u &=& -Df(x)^{-1}f(x),\\
N_f(x) &=& x \cos \| u \| + u \frac{\sin \| u \|}{\| u \|}.
\end{eqnarray*}
The size of the ball in Theorem \ref{th3} is equal to
$$R(f,\ze) = \min \left( \pi, \frac{3-\sqrt{7}}{2\ga(f, \ze)}\right).$$

\noindent {\bf Second example: the orthogonal group.} $\O_n$ denotes the orthogonal group. The tangent space at the identity matrix $id_n$ is equal to ${\cal A}_n$, the space of $n$ by $n$ antisymmetric matrices. More generally, the tangent space at $u \in \O_n$ is equal to
$$T_u \O_n = u {\cal A}_n.$$
This Riemannian structure is given by the usual scalar product of $n$ by $n$ matrices
$$\langle a,b \rangle = \trace(b^Ta)$$
for any $u \in \O_n$ and $a, b \in T_u \O_n$. The norm associated with this scalar product is the Frobenius norm and it is denoted by $\| a \|_F$, while the usual spectral norm is denoted by $\| a \|.$ $\O_n$ is a Lie group and this metric stucture is bi-invariant. Thus, the constant appearing in Definition \ref{def-K} is $K_u = 1$. 

The exponential map at $u \in \O_n$ is given by the exponential of matrices:
$$\exp_u(a) = u \exp{(u^{-1}a)}$$
for any $a \in T_u \O_n$, with 
$$ \exp{c} = \sum_{k=0}^\infty \frac{c^k}{k!}.$$
The inverse of the exponential is the logarithm
$$\log (id_n + b) = \sum_{k=1}^\infty (-1)^{k+1}\frac{b^k}{k}$$
defined for any matrix $b$ with $\| b \| < 1.$ Thus, the inverse of the exponential map
$$\exp_u^{-1}(b) = u \log (u^{-1}b)$$
is defined for any $b \in T_u \O_n$ such that $\| id_n - u^{-1}b \| < 1$ which is satisfied if and only if $\| u-b \| < 1.$
Consequently, the radius of injectivity is 
${\bf r}_u = 1$. 
Newton's method is given by 
$$N_f(u) = u \exp{(-u^{-1}Df(u)^{-1}f(u))}.$$
The size of the ball in Theorem \ref{th3} is equal to
$$R(f,\ze) = \min \left( 1, \frac{3-\sqrt{7}}{2\ga(f, \ze)}\right).$$

\noindent {\bf Third example: real projective space $\P_n(\R)$.}
Real projective space may be constructed as the quotient of 
$S^n \subset \R^{n+1}$ by the equivalence relation $x \equiv -x$.
Therefore, it has positive sectional curvature and hence $K_{x}=1$.
The radius of injectivity of the exponential is $\pi/2$.

Newton's method on $\P_n(\R)$ may be constructed as in the unit sphere 
(First example).

The size of the ball in Theorem \ref{th3} is equal to
$$R(f,\ze) = \min \left( \pi/2, \frac{3-\sqrt{7}}{2\ga(f, \ze)}\right).$$

\noindent {\bf Fourth example: Hermitian manifolds}

Let $M$ be an analytic, Hermitian $n$-dimensional manifold with 
metric $\langle \cdot , \cdot \rangle_H$. In particular, $M$ is
also a $2n$-dimensional analytic, Riemannian manifold with 
metric $\langle \cdot , \cdot \rangle = 
\mathbf{Re} \left( \langle \cdot , \cdot \rangle_H \right)$.

If $f: M \rightarrow \mathbb C^n$ is analytic, we define a real analytic
function $f_{\mathbb R}: M \rightarrow \mathbb R^{2n}$ by 
$f_{\mathbb R} (z) = \mathbf{Re}(f(z)), \mathbf{Im}(f(z))$.

Let $Df(z): T_zM \rightarrow \mathbb C^n$ denote the complex derivative
of $f$, in coordinates $z_1, \cdots, z_n$. Then,
\[
Df_{\mathbb R} =
\left[ \begin{array}{cc} 1 & i \\ 1 & -i \\ \end{array} \right]^{-1}
\left[ \begin{array}{cc} Df(z) & 0 \\ 0 & \overline{Df(z)}\\ \end{array} \right]
\left[ \begin{array}{cc} 1 & i \\ 1 & -i \\ \end{array} \right]
\]

It follows that $Df_{\mathbb R}(z) \left[\begin{array}{c}u \\ v\end{array}
\right]
= f_{\mathbb R}(z)$
if and only if $Df(z) \cdot (u+iv) = f(z)$. Therefore, Newton's method in
an Hermitian manifold is also given by
\[
N_f(z) = \mathrm{exp}_z \left( -Df(z)^{-1}f(z) \right)
\]

By the same argument, the invariants $\beta(f,z) = \| Df(z)^{-1} f(z) \|$
and $\gamma(f,z) = \sup_{k \ge 2} \left\| Df(z)^{-1} 
\frac{D^kf(z)}{k!}\right\|^{1/{k-1}}$ are equal, respectively, to
the Riemannian invariants $\beta(f_{\mathbb R}, z)$ and 
$\gamma(f_{\mathbb R},z)$. 

Therefore, Theorems~\ref{th3} to~\ref{thV4} apply verbatim to Hermitian
manifolds and maps $M \rightarrow \mathbb C^n$, or to vector fields
on Hermitian manifolds.

\section{Alternative formulation of the R-$\gamma$-Theorem}
\label{alternative}
In this section we investigate a question posed by an anonymous referee about the R-$\gamma$-Theorem (Theorem \ref{th3}). Using another proof we state it independently of the invariant $K(\zeta)$ introduced in Definition \ref{def-K}. We only state the R-$\gamma$-Theorem for mappings. The theorem
for vector fields is analogous.

\vskip 3mm
\begin{theorem} \label{th3M} {\bf (R$-\ga-$theorem)} 
There are constants $\nu_0 = 0.069778332 \ldots $ and $t_0 = 0.075262346 \ldots $ such that the following
statement is true. Let $\MM_n$ be geodesically complete and let $f: \MM_n \rightarrow \R^n$ be analytic.
Suppose that $f(\ze) = 0$ and $Df(\ze)$ is an isomorphism.

Let 
\[
R(f,\ze) = \min \left(
t_0 {\mathbf r}_{\zeta}, 
\nu_0 / \ga(f,\ze)
\right)
\]
If $d(z, \ze) \leq R(f,\ze)$,
then the Newton sequence $z_k = N_f^{(k)} (z)$ 
is defined for all $k \geq 0$,
and
$$ d(z_k, \ze) \leq \sigma \left ( \frac{1}{2} \right )^{2^k - 1}\beta(f,z).$$
\end{theorem}

Which theorem is the best? Theorem \ref{th3} or Theorem \ref{th3M}? 

When $\MM_n$ has a non-negative sectional curvature then, according to Corollary \ref{corthe3}, Theorem \ref{th3} gives a better result than Theorem \ref{th3M}. 
More generally, when $K_{\zeta} < \frac{1}{2 \nu_0} + \nu_0 - 2 = 5.235326440 \ldots $, the 
expression 
$\left( K_{\zeta} + 2 - \sqrt{K_{\zeta}^2 + 4 K_{\zeta} + 2} \right)/2$ in the hypothesis of Theorem~\ref{th3} is smaller than the constant $\nu_0$. This means that, unless geodesics spread away by a factor larger than 5 in
the relevant neighborhood, Theorem~\ref{th3} is sharper than Theorem~\ref{th3M}.
Otherwise Theorem~\ref{th3M} may be more useful.

We notice that, even if the formulation of Theorem \ref{th3M} doesn't depend on $K(\zeta)$, both radius of (proved) quadratic convergence depend on the metric at $\zeta$ via $\gamma (f, \zeta)$ and consequently on the curvature at this point. This also proves that, like in the case of linear spaces, the main invariant which estimates the size 
of the quadratic attraction basin of a root is the invariant gamma.

\vskip 3mm {\bf Proof of Theorem~\ref{th3M}:}
Let $\nu_0$ be the smallest positive root of the equation
$$\nu_0 / \psi(\nu_0)^2 = \alpha_0.$$ 
Numerically, $\nu_0 = 0.069778332 \ldots $
Also, let 
$$t_0 = \frac{s_0}{s_0 + \frac{1-\nu_0}{\psi(\nu_0)}} = 0.075262346 \dots$$
We assume that $\ze$ is such that $f(\ze)=0$ and $Df(\ze)$ is
an isomorphism. Let $z_0$ be such that $d(\ze, z_0) \le
\nu_0 / \gamma(f,\zeta)$. Since $\nu_0 < 1 - \sqrt{2}/2$,
\[
\nu = d(\ze, z_0) \gamma(f,\zeta) < 1 - \frac{\sqrt{2}}{2}
\]
and by Lemma~\ref{lem8} we have 
$$\beta(f,z_0) \le \frac{1-\nu}{\psi(\nu)} d(z_0,\zeta) = \frac{1-\nu}{\psi(\nu)} \frac{\nu}{\gamma(f,\zeta)}$$
and
$$\gamma(f,z_0) \le \frac{\gamma(f,\ze)}{(1-\nu) \psi(\nu)}.$$
Therefore,
\[
\alpha(f,z_0) \le \frac{\nu}{\psi(\nu)^2} \le \alpha_0.
\]

In order to apply the R$-\alpha-$Theorem, we need to show that
$\beta(f,z_0) \le s_0 {\mathbf r}_{z_0}$. Let $0<t<1$ be real number such that $d(z_0,\ze) \le t {\mathbf r}_\ze.$ Like previously
$$\beta(f,z_0) \le \frac{1-\nu}{\psi(\nu)} d(z_0,\zeta) \le \frac{1-\nu}{\psi(\nu)} t {\mathbf r}_\ze.$$
By Lemma \ref{lemn1}
$${\mathbf r}_\ze \le {\mathbf r}_{z_0} + d(z_0,\ze) \le {\mathbf r}_{z_0} + t {\mathbf r}_{\ze}$$
so that
$${\mathbf r}_\ze \le \frac{{\mathbf r}_{z_0}}{1-t}$$
and
$$\beta(f,z_0) \le \frac{1-\nu}{\psi(\nu)} \frac{t {\mathbf r}_{z_0}}{1-t}.$$
This gives $\beta(f,z_0) \le s_0 {\mathbf r}_{z_0}$ as soon as 
$$\frac{1-\nu}{\psi(\nu)} \frac{t }{1-t} \le s_0$$
or, equivalently,
$$t \le \frac{s_0}{s_0 + \frac{1-\nu}{\psi(\nu)}}$$
which is given by
$$t \le \frac{s_0}{s_0 + \frac{1-\nu_0}{\psi(\nu_0)}} = t_0 \simeq 0.075262346.$$

We can now apply the R-$\alpha$-theorem (Theorem~\ref{th4}):

\[
d(z_{k+1}, z_k) \le \left( \frac{1}{2} \right) ^{2^k - 1} \beta(f,z_0).
\]
Hence,

\[
d(z_k, \zeta) \le 
\left( \frac{1}{2} \right) ^{2^k - 1} 
\sigma \beta(f,z_0)
\]
\qed

\section{Conclusions and suggestions for further research}

In this paper, we gave a generalization of $\alpha$-theory 
for Riemannian (and therefore, Hermitian) manifolds. This 
generalization is subtle, due to the influence of new intrinsic
factors, such as the radius of injectivity of the exponential
and the curvature.

We developped an intrinsic approach avoiding the use of local charts or isometric imbeddings.
Except in the case of submanifolds, such imbeddings are often artificial and they lead to 
high dimensional problems, roughly speaking $n^3$ for a dimension $n$ manifold according to Nash's Embedding Theorem.

Our next objective is to implement this method. It is clear from the examples we have in mind and from the work already done that we have to take into account the data structure describing the considered problem. See for example
Celledoni-Iserles \cite{Celledoni} for Lie group methods, Edelman-Arias-Smith \cite{Edelman-Arias-Smith} for examples of manifolds described by the action of a group on a set and Adler-Dedieu-Margulies-Martens-Shub for a product of special orthogonal groups. These three papers show three different ways to compute the exponential map associated with 
the considered manifold and therefore three different ways to implement Newton's method.

\section{Acknowledgements}

Gregorio Malajovich was partially supported by Brazilian CNPq 
(Conselho Nacional de Desenvolvimento Cient\'{\i}fico e Tecnol\'ogico)
grant 300925/00-0(NV) and by
Funda\c c{\~ao} Jos\'e Pel\'ucio Ferreira. Part of this work
was done while Gregorio Malajovich and Jean-Pierre Dedieu 
were visiting the City University of Hong Kong (G.M. was 
supported by CERG grant no.~9040393).

The authors wish to thank the anonymous referees for their
comments and insights, J. Grifone, J.-M. Morvan and J.-M. Schlenker for many valuable discussions
about this subject.

\end{document}